\documentclass[11pt]{article}

\usepackage{lineno,hyperref,amsmath,amsfonts,amsthm,dsfont}

\theoremstyle{plain}
\newtheorem{theorem}{Theorem}[section]

\theoremstyle{definition}

\theoremstyle{remark}

\newcommand{\covth}[2]{\mathbb{C}\text{\em ov}\big({#1},{#2}\big)}
\newcommand{\cov}[2]{\mathbb{C}\text{ov}\big({#1},{#2}\big)}
\newcommand{\Idn}{\mbox{\textrm{$1\mskip-1.6\thinmuskip$I}} }

\newcommand{\cF}{{\mathcal F}}

\newcommand{\GR}{{\mathds{R}}}

\newcommand{\bE}{{\mathbb{E}}}
\newcommand{\bP}{{\mathbb{P}}}

\begin{document}
	\title{A complement to the Chebyshev integral inequality}

	\author{Adam Jakubowski\\[3mm]
		Nicolaus Copernicus University, Toru\'n, Poland\\[2mm]
		E-mail: adjakubo@mat.umk.pl
	}

	\date{}
	
	\maketitle
	
	\begin{abstract}	We give necessary and sufficient conditions for the Chebyshev inequality to be an equality.
	\end{abstract}

	
	\noindent{\em MSClassification 2010:} 60E15; 26D15.
	
	\section{Introduction}
	
	The simplest form of the Chebyshev integral (or algebraic) inequality holds for arithmetic means: if $x_1 \leq x_2 \leq \ldots \leq x_n$ and $y_1 \leq y_2 \leq \ldots \leq y_n$ are real numbers, then
	\begin{equation}\label{eq:1} \frac{1}{n} \big( x_1 y_1 + x_2 y_2 + \ldots + x_n y_n\big) \geq  \frac{1}{n} \big( x_1 + x_2 + \ldots + x_n\big)\cdot \frac{1}{n} \big( y_1 +  y_2 + \ldots + y_n\big).
	\end{equation}
	The continuous counterpart reads as follows: if $f, g: [a,b] \to \GR$ are non-decreasing, then
	\begin{equation}\label{eq:2} \frac{1}{b-a} \int_a^b f(x) g(x)\,dx  \geq  \frac{1}{b-a} \int_a^b f(x)\,dx \cdot \frac{1}{b-a} \int_a^b g(x)\,dx.
	\end{equation}
	Clearly, both (\ref{eq:1}) and (\ref{eq:2}) are particular cases of the following probabilistic statement.
	
	\begin{theorem}\label{th:chebyshev}
		If $X$ is a random variable on $\big(\Omega, \cF, \bP\big)$ and $f, g: \GR \to \GR$ are {\em non-decreasing}, and such that $\covth{f(X)}{g(X)}$ exists, then $f(X)$ and $g(X)$ are positively correlated, i.e.
		\begin{equation}\label{eq:3}
		\covth{f(X)}{g(X)} = \bE f(X)g(X) - \bE f(X) \bE g(X) \geq 0.
		\end{equation}
	\end{theorem}
	Notice that (\ref{eq:3}) valid for all non-decreasing $f$ and $g$  means that a single random variable is {\em associated}, what in turn is a cornerstone in the proof of association of independent random variables (see \cite{EPW67}, also \cite{BuSh07}).
	
	Somewhat surprisingly, a relatively recent work \cite{NiPe10} states that the Chebyshev inequality is equivalent to the classic Jensen inequality. In fact, both inequalities are ``dual'' in a specific sense (see \cite[Section 1.8]{NiPe06}).
	
	It is also remarkable that relation (\ref{eq:3}) admits direct consequences in some economic considerations - see \cite{Simo95}, \cite{Wage06}.
	
	We refer to \cite[Chapter IX]{MPF93} for the detailed report on developments related to the inequality that is nowadays called Chebyshev's. 
	
	The purpose of this note is to prove the following complement to Theorem \ref{th:chebyshev}.
	
	\begin{theorem}\label{th:zero}
		In assumptions of Theorem \ref{th:chebyshev}, $\covth{f(X)}{g(X)} = 0$ if, and only if, either $f(X)$ or $g(X)$ is a.s constant.
	\end{theorem}
	
	Taking into account the long history of the Chebyshev inequality the lack of explanation when the inequality is not strict seems to be unlikely. But we were not able to find any published reference proving Theorem \ref{th:zero} in full generality.
	
	Of course the trivial case (\ref{eq:1}) was clear in early eighteen-eighties due to the observation by A. Korkine (see \cite[p. 242]{MPF93}, also \cite[pp. 43--44]{HLP52}).
	
	\cite[p. 77]{Mitr72} states that the equality holds in (\ref{eq:2}) only if $f$ or $g$ is constant almost everywhere,  but the proof is missing. Notice that (\ref{eq:2}) corresponds to $X$ uniformly distributed on $[a,b]$. In a slightly more general case, when 
	the law of $X$ is given by a strictly positive density on $[a,b]$, 
	\cite[p. 40, Theorem 10]{Mitr70} refers to the original Chebyshev's proof \cite[pp. 128--131, 157--169]{Cheb48}. But Chebyshev worked under stronger assumptions (differentiability of $f$ and $g$). Likewise, \cite[Chapter IX]{MPF93} provide several results (by Winckler, Pickard, \ldots) where the strict inequality occurs in (\ref{eq:2}), but all of them are related to stronger assumptions imposed on functions $f$ and $g$. The reader may verify that in \cite[Chapter IX]{MPF93} no attention is paid to results similar to our Theorem \ref{th:zero}.
	
	Therefore it is surprising to find in \cite{Wage06} a statement (Theorem 1) referring to \cite[p. 248]{MPF93} and asserting that under Steffensen's assumptions on $f$ and $g$  the equality holds in (\ref{eq:3})  if, and only if, $f$ or $g$ ``are constant almost everywhere" (i.e. with respect to the Lebesgue measure).
	As simple examples (and our Theorem \ref{th:zero}) show this is incorrect. Notice that \cite{Wage06} refers also to the original paper \cite{Stef25a}.\footnote{Following \cite[p. 287]{MPF93} Wagener provides wrong year of publication of this article - 1920. The correct is 1925.} This article (written in Danish, a slightly extended English version is given in \cite{Stef25b}) contains strict inequalities in discrete form (hence for Riemann sums) and claims strict inequalities in integral form without further justification. In any case it is related to (\ref{eq:2}) rather than to (\ref{eq:3}).
	
	Providing a complete proof of Theorem \ref{th:zero} should prevent such mistakes and inaccuracies in the future. 
	
	In fact Theorem \ref{th:zero} may admit direct applications on its own.
	
	For example, let us consider a {\em non-degenerate} distribution function $F$ with finite variance. Let 
	\[ F^{\leftarrow}(u) = \inf \{ s\,;\, F(s) \geq u\}\]
	be its left-continuous inverse.  Then both $F^{\leftarrow}(u)$ and $- F^{\leftarrow}(1 - u)$ are {\em nondecreasing} in $u$ and therefore, if $U$ is uniformly distributed 
	on $[0,1]$, we have 
	\[\cov{F^{\leftarrow}(U)}{- F^{\leftarrow}(1-U)} > 0, \]
	or, equivalently,
	\[\cov{F^{\leftarrow}(U)}{F^{\leftarrow}(1-U)} < 0. \]
	We have justified a method of variance reduction, as described in \cite[Section 2.1]{BFS87}, without invoking the results of \cite{Whit76}.
	
	\section{Proofs}
	
	The proof of (\ref{eq:3}) is immediate, if we observe that 
	\[  \cov{f(X)}{g(X)} = 	\frac{1}{2}\bE \big(f(X) - f(Y)\big)\big(g(X) - g(Y)\big),\]
	where $Y$ is independent of $X$ and $Y \sim X$,
	and that for each $\omega$
	\[ \big(f(X(\omega)) - f(Y(\omega))\big)\big(g(X(\omega)) - g(Y(\omega))\big) \geq 0.\]
	
	To prove Theorem \ref{th:zero} suppose that
	\begin{align*}
	0 =\cov{f(X)}{g(X)} &= (1/2)  \bE \big(f(X) - f(Y)\big)\big(g(X) - g(Y)\big) \\
	&= (1/2) \bE_X\Big( \bE_Y \big(f(X) - f(Y)\big)\big(g(X) - g(Y)\big) \Big) 
	\end{align*}
	It follows that for $\bP_X$-almost all $x$
	$\bE_Y \big(f(x) - f(Y)\big)\big(g(x) - g(Y)\big) = 0$.
	Since still  $\big(f(x) - f(Y(\omega))\big)\big(g(x) - g(Y(\omega))\big) \geq 0$, we get
	\begin{equation}\label{eq:iter}
	\bP \big(\{ f(Y) = f(x)\} \cup \{ g(Y) = g(x)\}) =1, \ \text{ for $\bP_X$-almost all $x$}.
	\end{equation}
	
	Let 
	\[ A_f = \{ x\,;\, \bP\big( f(Y) = f(x)\big) > 0\}, \quad A_g = \{ x\,;\, \bP\big( g(Y) = g(x)\big) > 0\}.\]
	If $x \in A_f$, $f(x)$ is an atom of the distribution of $f(X)$. 
	Hence there are {\em distinct} numbers $u_1, u_2, \ldots $ such that $\bP\big( f(Y) = u_i\big) > 0$ and 
	\[ A_f = \bigcup_{i} f^{-1}(\{u_i\}).\]
	In particular, $A_f$ is measurable. If $\bP\big( Y \in A_f^c\big) > 0$, then by (\ref{eq:iter}) there exists $x_0$ such that $\bP\big(f(Y) = f(x_0)\big) = 0$ and 
	\[ 1 = \bP \big(\{ f(Y) = f(x_0)\} \cup \{ g(Y) = g(x_0)\}) = \bP \big( \{ g(Y) = g(x_0)\}).\]
	This proves the theorem. So we may and do assume that 
	\[ \bP\big( f(Y) \in \{u_1,u_2,\ldots\}\big) = \sum_i \bP\big(f(Y) = u_i\big) = 1.\]
	By symmetry we may also assume that for some {\em distinct} numbers $v_1, v_2, \ldots$
	\[ \bP\big( g(Y) \in \{v_1,v_2,\ldots\}\big) = \sum_k \bP\big(g(Y) = v_k\big) = 1.\]
	We can write	
	\begin{align*}
	f(X) &= \sum_{i=1}^{\infty} u_i \Idn_{\{f(X) = u_i\}},\ \
	f(Y) = \sum_{i=1}^{\infty} u_i \Idn_{\{f(Y) = u_i\}},\ \\ f(X) - f(Y) &= \sum_{i\neq j} (u_i - u_j) \Idn_{\{f(X) = u_i\}\cap \{f(Y) = u_j\}},\end{align*}
	and similarly
	\[ g(X) - g(Y) = \sum_{k\neq l} (v_k - v_l) \Idn_{\{g(X) = v_k\}\cap \{g(Y) = v_l\}}.\]
	
	Let $\omega \in\{f(X) = u_i, f(Y) = u_j, g(X) = v_k, g(Y) = v_l\}$, $i\neq j, k\neq l$.  Notice that $u_i > u_j$ implies $X(\omega) > Y(\omega)$, hence $v_k > v_l$ (the monotonicity of $f$ gives us only $v_k \geq v_l$, but we know that $v_k\neq v_l$). Similarly $u_i < u_j$ implies $v_k < v_l$. Therefore we always have $(u_i - u_j)(v_k-v_l) > 0$. We also have
	\begin{align*}
	0 &= \bE \big( f(X) - f(Y)\big)\big( g(X) - g(Y)\big)\\
	&= \sum_{i\neq j} \sum_{k\neq l} (u_i - u_j)(v_k- v_l) \bP\big( f(X) = u_i, f(Y) = u_j, g(X) = v_k, g(Y) = v_l\big)\\
	&= \sum_{i\neq j} \sum_{k\neq l} (u_i - u_j)(v_k- v_l) \bP\big( f(X) = u_i,  g(X) = v_k\big) \bP\big( f(Y) = u_j, g(Y) = v_l\big)\\
	&= \sum_{i\neq j} \sum_{k\neq l} (u_i - u_j)(v_k- v_l) p_{i,k} p_{j,l}.
	\end{align*}
	It follows that $p_{i,k}p_{j,l} = 0$, if $i\neq j$ and $ k \neq l$. 
	
	We have both
	\[ 1 = \sum_{i}\sum_{k} \bP\big( f(X) = u_i,  g(X) = v_k\big) = \sum_{i}\sum_k p_{i,k}\]
	and (keeping in mind that  $\Idn_{\{i=j\} \cup \{k=l\}} = \Idn_{\{i=j\}} + \Idn_{ \{k=l\}} - \Idn_{\{i=j, k=l\}}$)
	\begin{align*}
	1 &= \Big( \sum_{i,k} p_{i,k} \Big)^2 =  \sum_{i} \sum_k \sum_{j}\sum_l p_{i,k} p_{j,l}  = \sum_{i} \sum_k \sum_{j}\sum_l p_{i,k} p_{j,l} \Idn_{\{i=j\} \cup \{k=l\}}\\
	&=  \sum_{i} \sum_k \sum_l p_{i,k} p_{i,l}  + \sum_{i} \sum_k \sum_{j} p_{i,k} p_{j,k} -   \sum_{i} \sum_k p_{i,k}^2 \\
	&= \sum_{i} \sum_{k} p_{i,k} \Big(\sum_l p_{i,l} + \sum_j p_{j,k} - p_{i,k}\Big) =   \sum_{i} \sum_{k} p_{i,k} D_{i,k}
	\end{align*}
	Obviously $D_{i,k} \leq 1$ and if $p_{i,k} > 0$ then it must be $D_{i,k} = 1$. 
	
	Fo some $i_0,k_0$ $p_{i_0,k_0} > 0$. Then the whole mass of the joint distribution of $\big(f(Y),g(Y)\big)$ must be concentrated on the ``cross'' defined as the support of  $D_{i_0,k_0}$.
	If some $p_{r,k_0} > 0$, $r \neq i_0$, then by the repeated reasoning the complete mass of the distribution must be concentrated on the intersection of the two ``crosses", i.e. on the vertical axis containing $p_{i_0,k_0}$, i.e.
	\[ 1 = \sum_s p_{s,k_0} = \sum_s \bP\big( f(Y) = u_s,  g(Y) = v_{k_0}\big) = \bP\big( g(X) = v_k\big).\]
	Similarly, if for some  $q \neq k_0$ we have $p_{i_0,q} > 0$, then  $1 = \bP\big( f(Y) = u_{i_0}\big)$.
	If $p_{r,k_0} = 0$, $r \neq i_0$, and $p_{i_0,q} = 0$, $q \neq k_0$, then $1 = p_{i_0,k_0} = \bP\big ( f(X) = u_{i_0}, g(X) = v_{k_0}\big)$.

	\section*{Acknowledgements}
	The author would like to thank Thomas Mikosch and Boualem Djehiche for their help in accessing some of the old papers quoted in this work.
	

\end{document}